\documentclass[letterpaper, 11pt, conference]{ieeeconf}      

\IEEEoverridecommandlockouts                              
\usepackage{graphicx} 
\usepackage{epsfig} 
\usepackage{times} 
\usepackage{amsmath} 
\usepackage{amssymb}  
\usepackage{cleveref}
\usepackage{url}

\newcommand{\mathbbm}[1]{\text{\usefont{U}{bbm}{m}{n}#1}}

\usepackage{cite}
\newtheorem{theorem}{Theorem}

\newtheorem{assumption}[theorem]{Assumption}
\usepackage{color}
\usepackage{mathrsfs}
\usepackage{multirow}
\usepackage{caption}
\usepackage{subcaption}
\usepackage{cite}
\newtheorem{properties}[theorem]{Properties}

\newtheorem{definition}[theorem]{Definition}

{ 
\newtheorem{remark}[theorem]{Remark}

}
\usepackage{wrapfig}
\usepackage{float}

\newcommand{\bi}{\begin{itemize}}
\newcommand{\ei}{\end{itemize}}
\newcommand{\bd}{\begin{displaymath}}
\newcommand{\ed}{\end{displaymath}}
\newcommand{\be}{\begin{eqnarray*}}
\newcommand{\ee}{\end{eqnarray*}}

\title{\LARGE \bf 
Data-Driven Approximation of Transfer Operators: Naturally Structured  Dynamic Mode Decomposition
}
\author{Bowen Huang and Umesh Vaidya\\
\thanks{Financial support from the National Science Foundation grant ECCS-1002053 and ECCS 1150405 is gratefully acknowledged. Bowen Huang and Umesh Vaidya are with the Department of Electrical \& Computer Engineering,
Iowa State University, Ames, IA 50011
{\tt\small ugvaidya@iastate.edu}}%
}

\begin{document}
\maketitle
\begin{abstract}
In this paper, we provide a new algorithm
for the finite dimensional approximation of the linear
transfer Koopman and Perron-Frobenius operator from
time series data. We argue that existing approach for
the finite dimensional approximation of these transfer
operators such as Dynamic Mode Decomposition (DMD) and Extended Dynamic Mode Decomposition (EDMD) do not capture two important properties of these operators, namely positivity and Markov property. The algorithm
we propose in this paper preserve these two properties.
We call the proposed algorithm as naturally structured
DMD since it retains the inherent properties of these
operators. Naturally structured DMD algorithm leads to a better approximation of the steady-state dynamics of the system regarding computing Koopman and Perron-
Frobenius operator eigenfunctions and eigenvalues. However preserving positivity properties is critical for capturing the real transient dynamics of the system. This
positivity of the transfer operators and it's finite dimensional
approximation also has an important implication
on the application of the transfer operator methods for controller and estimator design for nonlinear systems from time series data.

\end{abstract}

%

\section{Introduction}

Transfer operator based methods involving Perron-Frobenius and Koopman operators are successfully applied for analysis and design of nonlinear dynamical systems in \cite{Dellnitz_Junge,Mezic2000,froyland_extracting,Junge_Osinga,Mezic_comparison,Dellnitztransport,mezic2005spectral,Mehta_comparsion_cdc,Vaidya_TAC,Vaidya_CLM_journal,raghunathan2014optimal,susuki2011nonlinear,mezic_koopmanism,mezic_koopman_stability,surana_observer}. The basic idea behind these methods is to shift the focus from the state space where the system evolution is nonlinear to measure space or space of functions where the system evolution is linear. The linearity of the transfer operator framework offers several advantages for analysis and design problems involving nonlinear systems. More importantly, this framework allows us to carry over our intuition from linear systems to nonlinear systems. One of the challenges in the application of these methods is in the computation of finite dimensional approximation of these infinite dimensional operators. \cite{dellnitz2002set}, proposed set-oriented numerical methods for the finite dimensional approximation of P-F transfer operators using the knowledge of system model. The data-driven model-free methods for the finite dimensional approximation of Koopman operator are proposed in \cite{Mezic2000}. Dynamic Mode Decomposition  (DMD) (\cite{DMD_schmitt}) and Extended DMD (\cite{rowley2009spectral,EDMD_williams}) are two of the popular algorithms that are proposed for approximating the spectrum (eigenvalues and eigenfunctions) of Koopman operator. These algorithms rely on mapping the time-series data from state space into space of observables using some finitely many basis functions in dictionary set. Finite dimensional approximation of the Koopman operator is then obtained as a matrix that best describes the evolution of the finite basis functions. The finite dimensional Koopman matrix is obtained as a solution of the least squares optimization problem. However, the existing approximation algorithm involving DMD and EDMD does not preserve some of the important properties of the Koopman, transfer operator. In particular, the Koopman operator is a positive operator, i.e., any positive function is mapped to positive function by Koopman operator \cite{Lasota}. Similarly, P-F and Koopman operators are adjoint operators. Furthermore, the P-F operator is a special class of Markov operator \cite{Lasota}. In the recent work by \cite{klus2015numerical}, the adjoint property of these two operators was exploited to provide a data-driven approximation of both Koopman and P-F operators.
 The Markov property of the P-F operator combined with the adjoint nature of two operators has important implications on the finite dimensional approximation of these two operators.

In this paper, we propose a new algorithm for the finite-dimensional approximation of these two operators that explicitly accounts for the positivity and Markov property and ensure that these features are retained in the finite-dimensional approximation. We show that preserving these properties allows one to better approximate the steady-state dynamics as captured by the spectrum (eigenvalue and eigenfunctions) of these operators but is essential to obtain the actual transient behavior of the system. We call the new algorithm for the finite dimensional approximation of transfer operator while preserving the properties of its infinite dimensional counterpart as Naturally Structured Dynamic Mode Decomposition (NSDMD).
We show that the problem of finding the finite dimensional approximation of the Koopman operator using NSDMD is a least square optimization problem with constraints and is convex. Using the adjoint property between the two transfer operators, we also construct the finite dimensional approximation of the P-F transfer operator. The P-F transfer operator is used to compute the finite dimensional approximation of the eigenfunction with eigenvalue one of the P-F operator capturing the steady-state invariant dynamics of the system. Structure preserving the property of our proposed NSDMD algorithm makes this possible. Furthermore, DMD and EDMD algorithm does not lead to stable finite dimensional Koopman matrix since the largest eigenvalue of the Koopman matrix is not guaranteed to be one. Since the Koopman operator obtained using NSDMD preserves the Markov property the largest eigenvalue is always one leading to a stable finite-dimensional approximation.

The organization of the paper is as follows. In Section \ref{section_operator}, we provide a brief overview of infinite dimensional operators and discuss the properties of these two operators. In Section \ref{section_numerics}, we present the overview of existing algorithms for the finite dimensional approximation of the operators, namely set-oriented numerics, DMD, and EDMD. We present the main result of this paper in the form of novel NSDMD algorithm that preserves the properties of these operators from infinite dimension to finite dimension in Section \ref{section_NSDMD}. Simulation results are presented in Section \ref{section_simulation} followed by conclusions in Section \ref{section_conclusion}.

\section{Transfer operators and their Spectrum}\label{section_operator}
Consider a discrete time dynamical system
\begin{eqnarray}
x_{t+1}=T(x_t)\label{system}
\end{eqnarray}
where $T:X\subset \mathbb{R}^N\to  X$ is assumed to be invertible and smooth diffeomorphism. Furthermore, we denote by ${\cal B}(X)$ the Borel-$\sigma$ algebra on $X$ and ${\cal M}(X)$ vector space of bounded complex valued measure on $X$.  Associated with this discrete time dynamical systems are two linear operators namely Koopman and Perron-Frobenius (P-F) operator. These two operators are defined as follows.
\begin{definition}[Perron-Frobenius Operator]  $\mathbb{P}_T:{\cal M}(X)\to {\cal M}(X)$ is given by
\[[\mathbb{P}\mu](A)=\int_{{\cal X} }\delta_{T(x)}(A)d\mu(x)=\mu(T^{-1}(A))\]
$\delta_{T(x)}(A)$ is stochastic transition function which measure the probability that point $x$ will reach the set $A$ in one time step under the system mapping $T$. \end{definition}
\begin{definition}[Invariant measures] are the fixed points of
the P-F operator $\mathbb{P}_T$ that are additionally probability measures. Let $\bar \mu$ be the invariant measure then, $\bar \mu$ satisfies
\[\mathbb{P}\bar \mu=\bar \mu\]
\end{definition}
Under the assumption that the state space $X$ is compact it is known that the P-F operator admits at least one invariant measure.
\begin{definition} [Koopman Operator] Given any $h\in\cal{F}$, $\mathbb{U}:{\cal F}\to {\cal F}$ is defined by
\[[\mathbb{U} h](x)=h(T(x))\]
\end{definition}

\begin{properties}\label{property}
Following properties for the Koopman and Perron-Frobenius operators can be stated.

\begin{enumerate}
\item [a).] For ${\cal F}=L_2(X,{\cal B}, \bar \mu)$ as the Hilbert space it is easy to see that
\begin{eqnarray*}
&&\parallel \mathbb{U}h\parallel^2=\int_X |h(T(x))|^2d\bar \mu(x)
\nonumber\\&=&\int_X | h(x)|^2 d\bar\mu(x)=\parallel h\parallel^2
\end{eqnarray*}
where we used the fact the $\bar \mu$ is an invariant measure. This implies that Koopman operator is unitary.

\item [b).] For any $h\geq 0$, we have $[\mathbb{U}h](x)\geq 0$ and hence Koopman is a positive operator.

\item [c).]For invertible system $T$, the P-F operator for the inverse system $T^{-1}:X\to X$ is given by $\mathbb{P}^*$ and $\mathbb{P}^*\mathbb{P}=\mathbb{P}\mathbb{P}^*=I$. Hence, the P-F operator is unitary.

\item [d).] If we define P-F operator act on the space of densities i.e., $L_1(X)$ and Koopman operator on space of $L_\infty(X)$ functions, then it can be shown that the P-F and Koopman operators are dual to each others as follows \footnote{with some abuse of notation we are using the same notation for the P-F operator defined on the space of measure and densities.}
\begin{eqnarray*}
&&\left<\mathbb{U} f,g\right>=\int_X [\mathbb{U} f](x)g(x)dx\nonumber\\&=&\int_Xf(y)g(T^{-1}(y))\left|\frac{dT^{-1}}{dy}\right|dy=\left<f,\mathbb{P} g\right>
\end{eqnarray*}
where $f\in L_{\infty}(X)$ and $g\in L_1(X)$ and the P-F operator on the space of densities $L_1(X)$ is defined as follows
\[[\mathbb{P}g](x)=g(T^{-1}(x))|\frac{dT^{-1}(x)}{dx}|\]

\item [e).] For $g(x)\geq 0$, $[\mathbb{P}g](x)\geq 0$.

Let $(X,{\cal B},\mu)$ be the measure space where $\mu$ is a positive but not necessarily the invariant measure of $T:X\to X$, then the P-F operator $\mathbb{P}:L_1(X,{\cal B},\mu)\to L_1(X,{\cal B},\mu)$  satisfies  following properties.

\item [f).] \[\int_X [\mathbb{P}g](x)d\mu(x)=\int_X g(x)d\mu(x)\]\label{Markov_property}
\end{enumerate}
\end{properties}

The linearity of the P-F operator combined with the properties \ref{property} (e) and \ref{property} (f), makes the P-F operator a particular case of Markov operator. This Markov property of P-F operator has significant consequences on its finite dimensional approximation. We will discuss this in section \ref{section_numerics} on set-oriented numerical methods for finite dimensional approximation of P-F operator.
Since $\mathbb{P}$ and $\mathbb{U}$ are unitary operators their spectrum lies on the unit circle. Given the adjoint nature of two operators, the spectrum of these operators are related. To study the connection between the spectrum of these two operators, we refer the interested readers to \cite{Mezic_comparison} and \cite{Mehta_comparsion_cdc} (Theorem 5 and Corollary 6) for results connecting the spectrum of transfer Koopman and P-F operator both in infinite dimensional and finite dimensional setting.

\section{Set-oriented numerics and Dynamic mode decomposition}\label{section_numerics}

\subsection{Set-oriented numerical methods}\label{section_setoriented}
Set oriented numerical methods are primarily developed for the finite dimensional approximation of the Perron-Frobenius operator for the case where system dynamics are known \cite{dellnitz2002set, GAIO01}. However, these algorithms can be modified or extended to the case where system information is available in the form of time series data. The basic idea behind set-oriented numerics is to partition the state space, $X$, into the disjoint set of boxes $D_i$ such that $X=\cup_{i=1}^N D_i$. Consider a finite partition $X^{'}=\{D_1,\ldots, D_K\}$.
Now, instead of a Borel $\sigma$-algebra, consider a $\sigma$-algebra of all possible subsets of $X$. A real-valued measure $\mu_j$ is defined by ascribing to each element $D_j$ a real number. This
allows one to identify the associated measure space with a
finite-dimensional real vector space $\mathbb{R}^K$. A given mapping
$T : X \to X$ defines a stochastic transition function $\delta_{T(x)}(\cdot)$. This function can be used to obtain a coarser representation of P-F operator denoted by ${\bf P}':\mathbb{R}^{K\times K}\to \mathbb{R}^{K\times K}$ as follows: For $\mu^{'}=(\mu_1^{'},\ldots, \mu_K^{'})$ we define a measure on $X$ as
\[d\mu(x)=\sum_{k=1}^K \mu_k^{'}\chi_{D_k}(x)\frac{dm(x)}{m(D_k)}\]
where $\chi_{D_k}(x)$ is the indicator function of $D_k$ and $dm$ is the Lebesgue measure. The finite dimensional approximation of the P-F matrix, ${\bf P}'$, can now be obtained as follows:
\begin{eqnarray}&&\nu_i'=[{\bf P}'\mu'](D_i)=\sum_{j=1}^K \int_{D_j}\delta_{T(x)}(D_i)\mu_j'\frac{dm(x)}{m(D_j)}
\nonumber\\&=&\sum_{j=1}^K \mu_k'{\bf P}'_{ij}\end{eqnarray}
where
\[{\bf P}'_{ij}=\frac{m(T^{-1}(D_j)\cap D_i)}{m(D_j)}\]
The resulting matrix $\bf {P}'$ is a Markov matrix and is row stochastic if we consider state $\mu'$ to be a row vector multiplying from the left of $P$. The individual entries of this Markov matrix can be obtained by  Monte-Carlo approach by running simulation over short time interval starting from different initial conditions. Typically individual boxes $D_i$ will be populated with $M$ uniformly distributed initial conditions. The entry ${\bf P}_{ij}$ is then approximated by fraction of initial conditions that are in box $D_j$ in one forward iteration of the mapping $T$. The Monte Carlo based approach can be extended for computation of the P-F transfer operator from time series data. Let $\{x_0,T(x_0),\ldots, T^{N-1}(x_0)\}$ be the time series data set. The number of initial conditions in box $i$ is then given by
\[\sum_{k=0}^{N-1} \chi_i(T^k(x_0)) \]
where $\chi_i$ is the indicator function of box $i$. The $(i,j)$ entry for P-F matrix ${\bf P}'_{ij}$ is then given by the fraction of these initial conditions from box $i$ that ends up in box $j$ after one iterate of time and is given by following formula.
\[{\bf P}'_{ij}=\frac{1}{\sum_{k=0}^{N-1} \chi_i (T^{k}(x_0))}\sum_{k=0}^{N-1}\chi_i(T^{k}(x_0))\chi_j(T^{k+1}(x_0)).\]

\subsection{Dynamic mode decomposition (DMD) and Extended DMD}
Dynamic Mode Decomposition method (DMD) has been introduced \cite{DMD_schmitt} for the dynamical analysis of the fluid flow field data. In the context of this paper, DMD can be viewed as a computation algorithm for approximating the spectrum of Koopman operator \cite{rowley2009spectral}. Extension of the DMD is presented in the form of Extended DMD (EDMD) which does a better job in approximating the spectrum of Koopman operator for both linear and nonlinear underlying system. In the following, we briefly explain the EDMD algorithm and show how the solution of DMD algorithm can be derived as a special case of EDMD. Consider snapshots of data set obtained from simulating a discrete time dynamical system $z\to T(z)$ or from an experiment
\begin{eqnarray}
\overline X = [x_1,x_2,\ldots,x_M],&\overline Y = [y_1,y_2,\ldots,y_M] \label{data}
\end{eqnarray}
where $x_i\in X$ and $y_i\in X$. The two pair of data sets are assumed to be two consecutive snapshots i.e., $y_i=T(x_i)$. Now let $\mathcal{D}=
\{\psi_1,\psi_2,\ldots,\psi_K\}$ be the set of dictionary functions or observables. The dictionary functions are assumed to belong to $\psi_i\in L_2(X,{\cal B},\mu)={\cal G}$, where $\mu$ is some positive measure not necessarily the invariant measure of $T$. Let ${\cal G}_{\cal D}$ denote the span of ${\cal D}$ such that ${\cal G}_{\cal D}\subset {\cal G}$. The choice of dictionary functions are very crucial and it should be rich enough to approximate the leading eigenfunctions of Koopman operator. Define vector valued function $\mathbf{\Psi}:X\to \mathbb{C}^{K}$
\begin{equation}
\mathbf{\Psi}(\boldsymbol{x}):=\begin{bmatrix}\psi_1(x) & \psi_2(x) & \cdots & \psi_K(x)\end{bmatrix}
\end{equation}
In this application, $\mathbf{\Psi}$ is the mapping from physical space to feature space. Any function $\phi,\hat{\phi}\in \mathcal{G}_{\cal D}$ can be written as
\begin{eqnarray}
\phi = \sum_{k=1}^K a_k\psi_k=\boldsymbol{\Psi^T a},\quad \hat{\phi} = \sum_{k=1}^K \hat{a}_k\psi_k=\boldsymbol{\Psi^T \hat{a}}
\end{eqnarray}
for some set of coefficients $\boldsymbol{a},\boldsymbol{\hat{a}}\in \mathbb{C}^K$. Let \[ \hat{\phi}(x)=[\mathbb{U}\phi](x)+r,\]
where $r\in\mathcal{G}$ is a residual function that appears because $\mathcal{G}_{\cal D}$ is not necessarily invariant to the action of the Koopman operator. To find the optimal mapping which can minimize this residual, let $\bf K$ be the finite dimensional approximation of the Koopman operator. Then the  matrix $\bf K$ is obtained as a solution of least square problem as follows 
\begin{equation}\label{edmd_op}
\min\limits_{\bf K}\parallel {\bf G}{\bf K}-{\bf A}\parallel_F
\end{equation}
\begin{eqnarray}\label{edmd1}
&&{\bf G}=\frac{1}{M}\sum_{m=1}^M \boldsymbol{\Psi}({x}_m)^\top \boldsymbol{\Psi}({x}_m)\nonumber\\
&&{\bf A}=\frac{1}{M}\sum_{m=1}^M \boldsymbol{\Psi}({x}_m)^\top \boldsymbol{\Psi}({y}_m),
\end{eqnarray}
with ${\bf K},{\bf G},{\bf A}\in\mathbb{C}^{K\times K}$. The optimization problem (\ref{edmd_op}) can be solved explicitly to obtain following solution for the matrix $\bf K$
\begin{eqnarray}
{\bf K}_{EDMD}={\bf G}^\dagger {\bf A}\label{EDMD_formula}
\end{eqnarray}
where ${\bf G}^{\dagger}$ is the psedoinverse of matrix $\bf G$.
Hence, under the assumption that the leading Koopman eigenfunctions are nearly contained within $\mathcal{G}_{\mathcal{D}}$, the subspace spanned by the elements of $\mathcal{D}$. The eigenvalues of $\bf K$ are the EDMD approximation of Koopman eigenvalues. The right eigenvectors of $\bf K$ generate the approximation of the eigenfunctions in (\ref{EDMD_eigfunc_formula}). In particular, the approximation of Koopman eigenfunction is given by
\begin{equation}\label{EDMD_eigfunc_formula}
\phi_j=\boldsymbol{\Psi} v_j
\end{equation}
where $v_j$ is the $j$-th right eigenvector of $\bf K$, $\phi_j$ is the eigenfunction approximation of Koopman operator associated with j-th eigenvalue.

DMD is a particular case of EDMD, and it corresponds to the case where the dictionary functions are chosen to be equal to ${\cal D}=\{e_1^\top,\ldots, e_K^\top\}$, where $e_i\in \mathbb{R}^N$ is a unit vector with $1$ at $i^{th}$ position and zero elsewhere. With this choice of dictionary function, it can be shown the approximation of the Koopman operator using DMD approach can be written as
\[{\bf K}_{DMD}=\overline Y\;\overline X^{\dagger},\]
where $\overline X$ and $\overline Y$ are data set as defined in (\ref{data}).
\section{Naturally structured dynamic mode decomposition}\label{section_NSDMD}
In this section, we provide a new algorithm for the finite dimensional approximation of the Koopman and P-F operator that preserves some of the properties of these two operators. In particular, we develop an algorithm that preserves the positivity property of the Koopman operator. Furthermore, the adjoint nature of Koopman and P-F operators is used to impose additional constraints on the entries of the Koopman operator. These structural properties are not considered in the existing algorithms involving DMD and EDMD for the finite dimensional approximation of the Koopman operator.
We show using examples that preserving these properties leads to a better approximation of eigenfunctions and eigenvalues of the transfer operators, but these features are essential to capture the correct transient behavior of the system. Capturing real transient dynamics is of particular importance towards the applications of the transfer operator for data-driven control and estimation problems.

In our proposed numerical algorithm for finite dimensional approximation of transfer operators from data we start with the choice of dictionary functions ${\cal D}=\{\psi_1,\ldots,\psi_K\}$, where $\psi_i(x)\in {\cal G}=L_2(X,{\cal B},\mu)$. As already stated the choice of dictionary function is crucial and should be rich enough to approximate the Koopman eigenfunctions. Similarly the data set generated by the dynamics should be rich enough to carry the information about the inherent dynamics of the system. We believe that the proper choice of dictionary function and data set are intimately connected.

We made following assumptions on the choice of dictionary function.
\begin{assumption}\label{assumption_dic}
We assume that the dictionary function $\psi_i(x)\geq 0$ for $i=1,\ldots, K$ and the inner product $\Lambda$ of the dictionary functions, $\Lambda=\langle\boldsymbol{\Psi}(x),\boldsymbol{\Psi}(x)\rangle$ with $[\Lambda]_{ij}=\langle\psi_i,\psi_j\rangle$ is symmetric positive definite  matrix.
\end{assumption}
\begin{remark} Gaussian radial basis function (RBF) given by $\exp^{-\frac{\parallel x-x_i \parallel}{\sigma^2}}$, serves as a good approximation for the choice of dictionary functions satisfying the above assumption.
\end{remark}
Let ${\cal G}_{\cal D}$ be the span of these dictionary functions. Now consider any function  $\phi$ and $\hat \phi$ in ${\cal G}_{\cal D}$, we can express these functions as
\begin{eqnarray}
\phi = \sum_{k=1}^K a_k\psi_k=\boldsymbol{\Psi^T a},\quad \hat{\phi} = \sum_{k=1}^K \hat{a}_k\psi_k=\boldsymbol{\Psi^T \hat{a}}
\end{eqnarray}
Again function $\phi$ and $\hat \phi$ are related as follows
\[\hat \phi(x)= [\mathbb{U}\phi](x)+r\]
where $r\in {\cal G}$ and represents the error and arise because of the fact that ${\cal G}_{\cal D}$ is not necessarily invariant under the action of Koopman operator. The extended DMD seeks to find the matrix ${\bf K}\in \mathbb{R}^{K\times K}$ that does the best job in mapping $\boldsymbol{a}$ to $\boldsymbol{\hat a}$. The matrix $\bf K$ is obtained as a solution of the least square problem as outlined in Eqs. (\ref{EDMD_formula}) and (\ref{edmd1}).  Now consider a case where $\phi(x)\geq 0$. Then under Assumption \ref{assumption_dic}, we know that $a_i\geq 0$. Using the positivity property of the Koopman operator, we know that $[\mathbb{U}\phi](x)\geq 0$. The vector $\boldsymbol{a}$ is mapped to $\hat {\boldsymbol a}$ by the finite dimensional matrix $\bf K$. To preserve the positivity property of the Koopman operator (i.e., property \ref{property}b) we require that coefficient $\hat {a}_i$ are also positive. This, in turn, implies that the mapping $\bf K$ should satisfy the property

\begin{eqnarray}{\bf K}_{ij}\geq 0,\;\;{\rm for}\;\; i,j=1,\ldots, K. \label{positive}
\end{eqnarray}
Let $\bf P$ be the finite dimensional approximation of the P-F operator. Since P-F is Markov operator, its finite dimensional approximation constructed on the dictionary function satisfying Assumption \ref{assumption_dic} has some properties. In particular, consider any density function, $\varphi$, expressed as linear combinations of dictionary functions
\[\varphi=\sum_{k=1}^K \boldsymbol{b}_k \psi_k,\;\;\;\boldsymbol{b}_k\geq 0.\]
We have
\[[\mathbb {P}\varphi](x)=\hat \varphi(x)+r=\sum_{k=1}^K \boldsymbol{\hat{b}}_k \psi_k+r,\]
where $r\in {\cal G}$ is the residual term which arise because ${\cal G}_{\cal D}$ is not invariant under the action of the P-F operator. The finite dimensional approximation of the P-F operator, $\bf P$ maps coefficient vector $\boldsymbol{b}$ to $\boldsymbol{\hat{b}}$, i.e., $\boldsymbol{\hat{b}}={\bf P}\boldsymbol{b}$

We are interested in approximating P-F operator such that the Markov property \ref{property}(f) of the infinite dimensional P-F operator is preserved. Since $[\mathbb {P}\varphi](x)\geq 0$ we have $\boldsymbol{b}_k\geq 0$ for all $k$. Hence for preserving the Markov property we require that
\begin{eqnarray}\boldsymbol{b}^\top {\bf 1}=\boldsymbol{\hat{b}}^\top {\bf1},\label{markov}\end{eqnarray}
where $\bf 1$ is a vector of all ones. 

Based on the adjoint property of Koopman and P-F operators, 
we have
\[\langle\mathbb{U}\phi,\varphi\rangle=\langle\phi,\mathbb{P}\varphi\rangle\]
Writing $\varphi$ and $\phi$ as linear combinations of basis function and using the definition of inner product from Assumption \ref{assumption_dic}, we can approximate the adjoint relationship as follows:
\begin{eqnarray}
\nonumber\langle\mathbb{U}\phi,\varphi\rangle\cong (\bf{K}\boldsymbol{a})^\top\Lambda\boldsymbol{b}&,& \langle\phi,\mathbb{P}\varphi\rangle\cong \boldsymbol{a}^\top\Lambda{\bf P}\boldsymbol{b}\\
\boldsymbol{a}^\top\bf{K}^\top\Lambda\boldsymbol{b}&=&\boldsymbol{a}^\top\Lambda\bf{P}\boldsymbol{b}
\end{eqnarray}
Since above is true for all $\boldsymbol{a}$ and $\boldsymbol{b}$, we have $\bf{K}^\top\Lambda=\Lambda\bf{P}$. 
Combining (\ref{positive}), (\ref{markov}) and the adjoint property of P-F and Koopman operator (i.e., ${\bf P}^\top=\Lambda{\bf K}\Lambda^{-1}$), it follows that for the finite dimensional approximation of the transfer operator to preserve the positivity and Markov properties of its infinite dimensional counterpart then $\bf K$ should satisfy following conditions
 \[{[\Lambda{\bf K}\Lambda^{-1}]}_{ij}\geq 0,\;\;\;\sum_{j=1}^K {[\Lambda{\bf K}\Lambda^{-1}]}_{ij}=1,\;i,j=1,\ldots, K.\]
 This leads to the following optimization based formulation for the computation of matrix $\bf K$
\begin{eqnarray}\label{optimization_problem}
\min\limits_{\bf K} & \parallel {\bf G}{\bf K}-{\bf A}\parallel_F\\\nonumber
\text{subject to} & {\bf K}_{ij} \geq 0\\\nonumber
& [{\Lambda {\bf K}\Lambda^{-1}}]_{ij}\geq 0\\\nonumber
& \Lambda{\bf K}\Lambda^{-1}\mathbbm{1} = \mathbbm{1}
\end{eqnarray}
where $\bf G$ and $\bf A$ are defined as follows:

\begin{eqnarray}\label{edmd2}
&&{\bf G}=\frac{1}{M}\sum_{m=1}^M \boldsymbol{\Psi}({x}_m)^\top \boldsymbol{\Psi}({x}_m)\nonumber\\
&&{\bf A}=\frac{1}{M}\sum_{m=1}^M \boldsymbol{\Psi}({x}_m)^\top \boldsymbol{\Psi}({y}_m),
\end{eqnarray}
with ${\bf K},{\bf G},{\bf A}\in\mathbb{C}^{K\times K}$ and the data set snapshots $\{x_n,y_n\}$ as defined in (\ref{data}). The optimization problem (\ref{optimization_problem}) is a convex and can be solved using one of the standard optimization toolbox for solving convex problem.

It is important to emphasize that the matrix $\bf K$ serves two purposes; a) approximation of Koopman operator if we multiply vector from right; b) approximation to P-F operator if we multiply vector from left.
\[{\rm Koopman \;operator} \;\;v_{t+1}={\bf K}v_t\]
\[{\rm P-F \;operator}\;\;u_{t+1}=u_t\bf P\]
where ${\bf P} =\Lambda{\bf K}\Lambda^{-1}$, $v_t\in \mathbb{R}^K$ is column vector and $u_t\in \mathbb{R}^K$ is row vector, and $t$ is the time index.

Since $\bf P$ is row stochastic, it is guaranteed to have at least one eigenvalue one. Let, $\bar u_1$ be the left eigenvector with eigenvalue one  of the $\bf P$ matrix. Then the approximation to the invariant density for the dynamical system, $T$, i.e., $\bar \varphi_1(x)$, can be obtained using following formula
\[\bar \varphi_1(x)=\boldsymbol\Psi(x)\bar u_1^\top.\]
Eigenfunction with eigenvalue $\lambda$  can be obtained as $\bar \varphi_{\lambda}=\boldsymbol\Psi(x)\bar u_{\lambda}^\top$, where $\bar u_{\lambda}^\top$ is the left eigenvector with eigenvalue $\lambda$ of matrix $\bf P$. Koopman eigenfunction with eigenvalue $\lambda$. We will refer to these eigenfunctions obtained using the left eigenvector of the $\bf P$ matrix as P-F eigenfunction. Similarly, approximate eigenfunctions of Koopman operator can be obtained using the right eigenvector of the $\bf K$ matrix. Let $\bar v_\lambda$ be the right eigenvector with eigenvalue $\lambda$ of the $\bf K$ matrix then the approximate Koopman eigenfunction $\bar \vartheta_\lambda$ can be obtained as  follows:
\[\bar \vartheta_\lambda(x)=\boldsymbol{\Psi}(x)\bar v_\lambda.\]

We show that NSDMD preserve the stability property of the original system and this is one of the main advantages of the proposed algorithm.  In particular, that certificate in the form of Lyapunov measure can be computed using the $\bf K$ matrix. \cite{Vaidya_TAC} introduced the Lyapunov measure for almost everywhere stability verification of general attractor set in the nonlinear dynamical system.  The Lyapunov measure is computed using transfer operator-based framework. \cite{Vaidya_TAC} utilized set-oriented numerical methods for the finite dimensional approximation of the P-F operator from system dynamics. However, data-driven approach for verifying the stability of attractor set will involve making use of matrix $\bf K$ for computing Lyapunov measure. The procedure for calculating the Lyapunov measure will remain the same; the only change is that instead of using the P-F matrix constructed using set-oriented numerical method one can use the $\bf K$ build from time series data. In the simulation section, we present results for the computation of stability certificate. Different optimization problems can be formulated based on the main optimization formulation in Eq. (\ref{optimization_problem}). These different optimization formulation will try to preserve one or all the properties of these two operators. In particular, we have following different cases.\\
\textbf{Case I}: With positivity constraint on ${\bf K}$ only
\begin{eqnarray}\label{optproblem1}
\min\limits_{\bf K} & \parallel {\bf G}{\bf K}-{\bf A}\parallel_F\\\nonumber
\text{subject to} & {\bf K}_{ij} \geq 0\\\nonumber
\end{eqnarray}

\textbf{Case II}: With positivity and Markov constraint on ${\bf P}$ only
\begin{eqnarray}\label{optproblem2}
\min\limits_{\bf K} & \parallel {\bf G}{\bf K}-{\bf A}\parallel_F\\\nonumber
\text{subject to} & [{\Lambda {\bf K}\Lambda^{-1}}]_{ij}\geq 0\\\nonumber
& \Lambda{\bf K}\Lambda^{-1}\mathbbm{1} = \mathbbm{1}
\end{eqnarray}
Both the optimization formulation (\ref{optproblem1}) and (\ref{optproblem2}) are convex formulation. 

{\bf Case III}: This case corresponds to combining both Case I and Case II and the optimization formulation corresponding to this case is given in Eq. (\ref{optimization_problem}).


\section{simulation results}\label{section_simulation}
The simulation results in this section are obtained by solving the optimization problems using  GUROBI solver coded in MATLAB.\\
{\underline{\it 2D system}}: For this example we use optimization formulation from {\bf Case I}. A simple 2D nonlinear system is considered first. The differential equation of the system is given as follows,
\begin{eqnarray}\nonumber
\dot x &=& x-x^3+y\\
\dot y &=& 2x-y\label{2d_system}
\end{eqnarray}

This continuous time system has 2 stable equilibrium points, located at $(\pm\sqrt{3},\pm2\sqrt{3})$ and one saddle point at $(0,0)$. To generate time-series data of $T=10$, $1000$ initial conditions from $[-5,5]\times[-5,5]$ are randomly chosen and propogated using ode23t solver in MATLAB, sampled by $\Delta t=0.1$. The naturally structured dynamic mode decomposition (NSDMD) algorithm is then implemented with Gurobi solver. The following simulation results are obtained with 500 dictionary functions and $\sigma=0.45$.

\begin{figure}[h!]
\centering
\includegraphics[width=.9\linewidth]{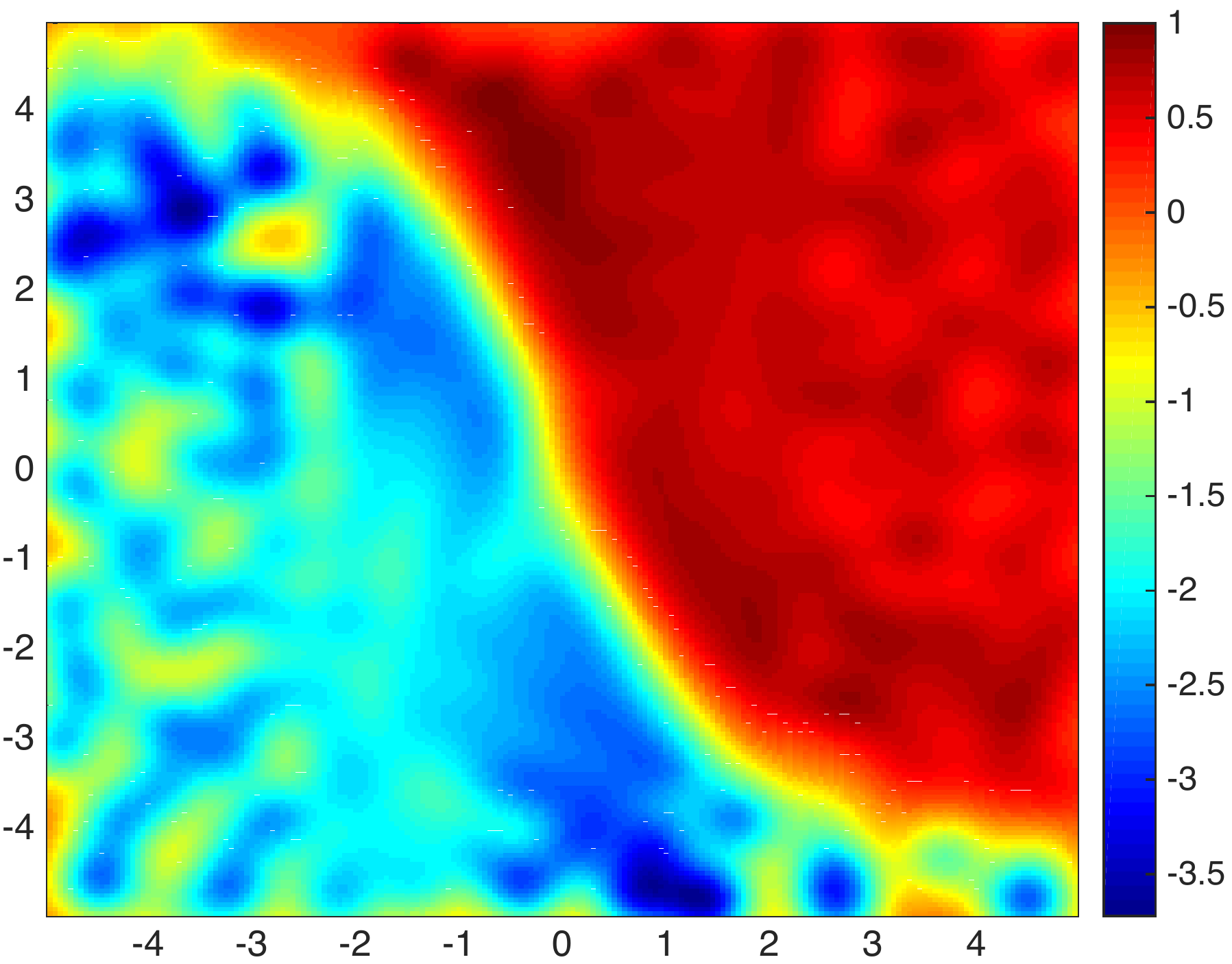}
\caption{\small {CASE-I: Koopman eigenfunction for eigenvalue $1$ for system (\ref{2d_system}) using NSDMD}}\label{2D_koopman1}
\end{figure}

\begin{figure}[h!]
\centering
\includegraphics[width=.9\linewidth]{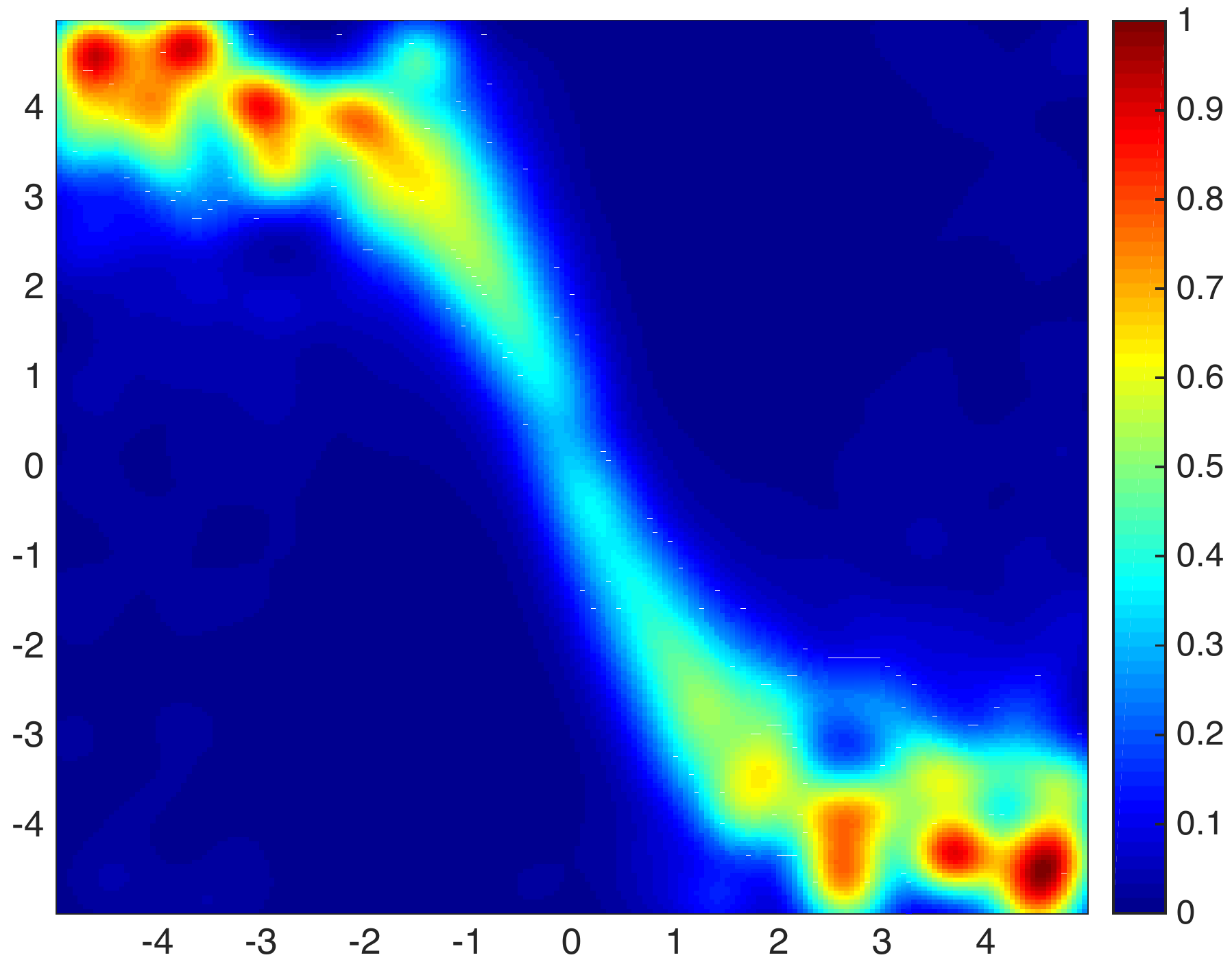}
\caption{\small {CASE-I: Koopman eigenfunction for eigenvalue $0.97$ for system (\ref{2d_system}) using NSDMD}}\label{2D_koopman2}
\end{figure}

In Fig. \ref{2D_koopman1} and Fig. \ref{2D_koopman2}, we plot the Koopman eigenfunctions associated with eigenvalue 1 using NSDMD algorithm. The eigenfunction with eigenvalue one is clearly shown to separate the two domain of attraction. The separatrix region separating the two domain of attractions is captured by the eigenfunction with second dominant eigenvalue.


{\underline{\it Duffing Oscillator}}: The simulation results for this example is obtained using formulation of {\bf Case I}. The duffing oscillator is given by following differential equation.
\begin{eqnarray}
\ddot x=-0.5 \dot x-(x^2-1)x
\end{eqnarray}

\begin{figure}[h!]
\centering
\includegraphics[width=.9\linewidth]{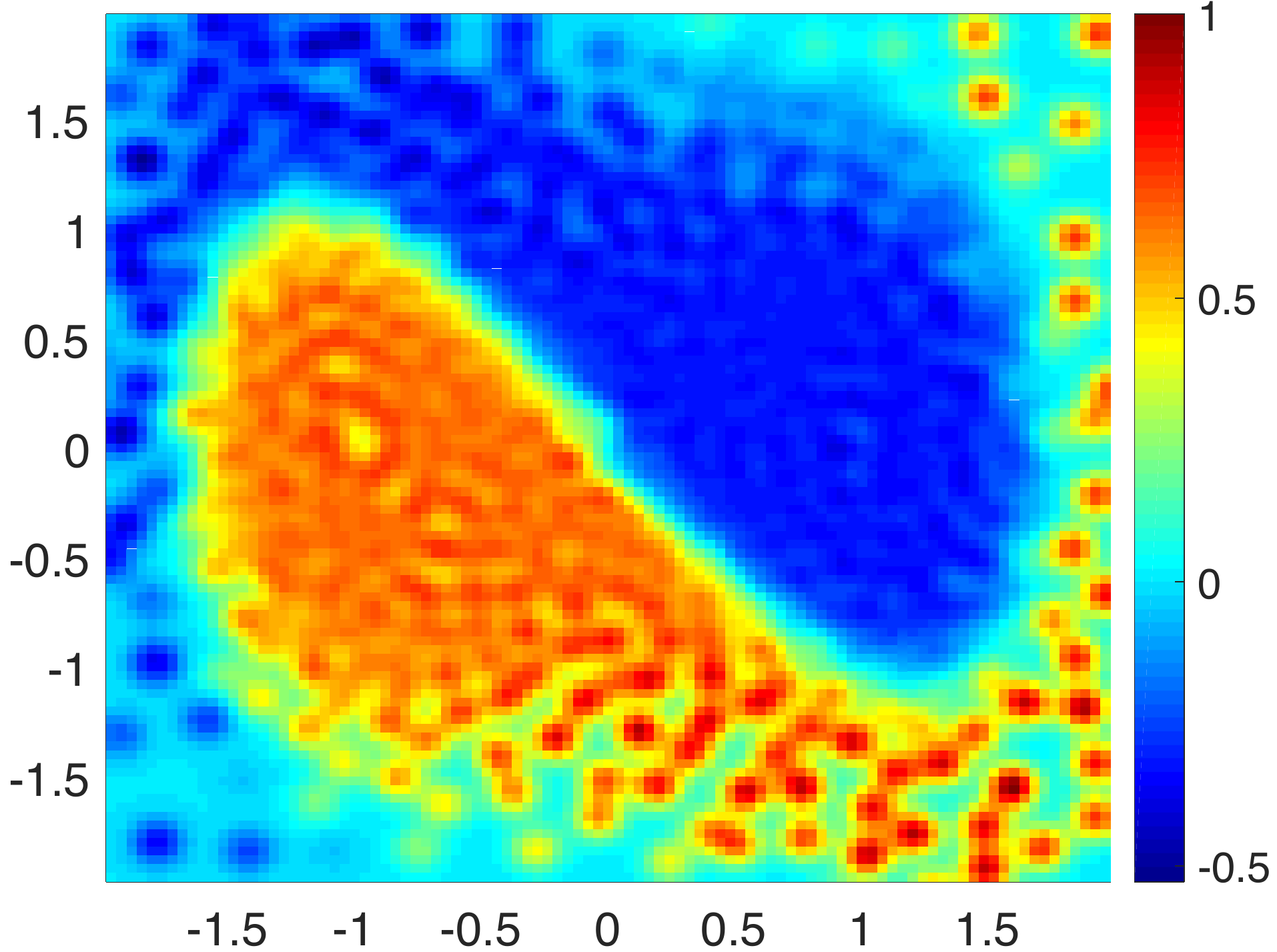}
\caption{\small {CASE-I: Koopman eigenfunction for eigenvalue $1$ for Duffing oscillator using NSDMD}}\label{duffing_koopman1}
\end{figure}

\begin{figure}[h!]
\centering
\includegraphics[width=.9\linewidth]{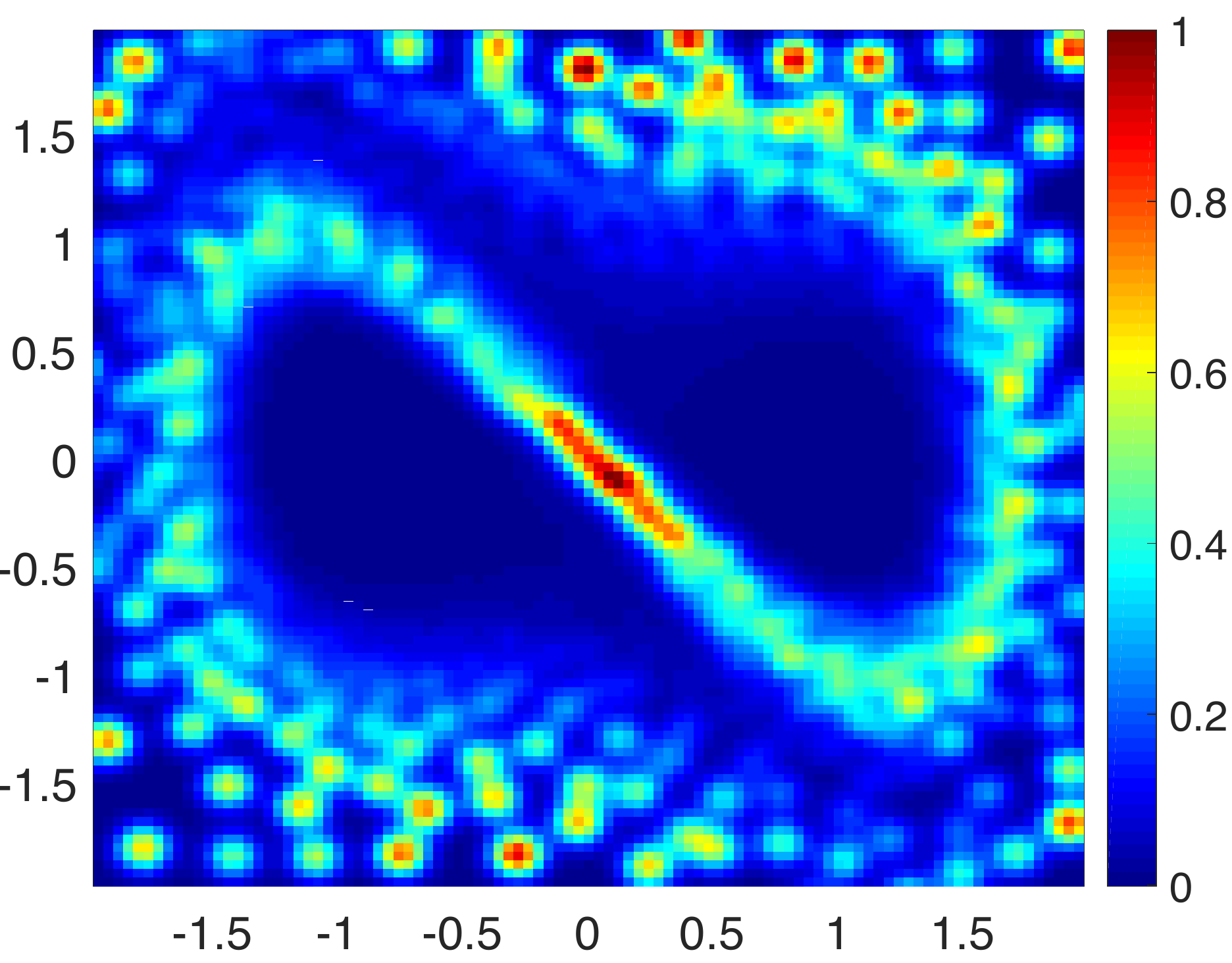}
\caption{\small {CASE-I: Koopman eigenfunction for eigenvalue $0.93$ for Duffing oscillator using NSDMD}}\label{duffing_koopman2}
\end{figure}

The time step for the continuous time system is chosen to be equal to $\Delta t=0.25$ with the total period of $T=2.5$ and $1000$ randomly chosen initial conditions. We solve the differential equation in MATLAB with $ode45$ solver. We use $500$ Gaussian radial basis functions to form the dictionary set with $\sigma=0.1$. In Fig. \ref{duffing_koopman1} and Fig. \ref{duffing_koopman2}, we plot the first two dominant eigenfunctions of the Koopman operator obtained using NSDMD algorithm. Similar to the example 1, we notice the first two dominant Koopman eigenfunctions carry information about the domain of attraction of the two equilibrium point.

\begin{figure}[h!]
\includegraphics[width=.9\linewidth]{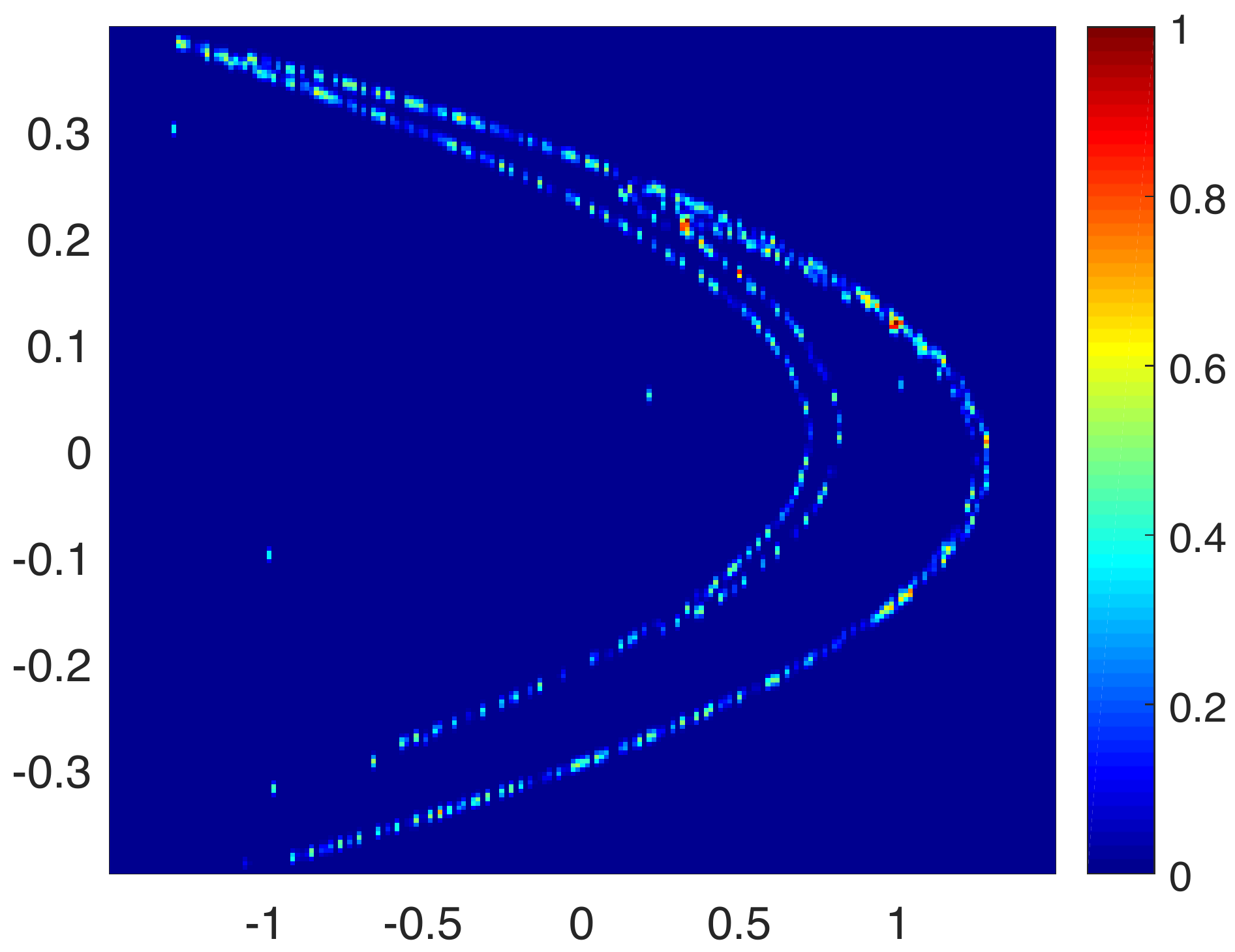}
\centering
\caption{\small {CASE-II: P-F eigenfunction for eigenvalue $1$ for Henon map using NSDMD}}
\label{henon_pf1}
\end{figure}

{\underline{\it Henon Map}}: Consider a following discrete-time system for the Henon map
\begin{eqnarray}
x_{t+1}&=&1-ax_t^2+y_t\nonumber\\
y_{t+1}&=&b x_t
\end{eqnarray}
with $a=1.4$ and $b=0.3$. Time series data starting from one initial condition over $5000$ time step is generated. Dictionary set is constructed using $500$ Gaussian radial basis functions. $K$-means clustering method is used for selecting the centers of these Gaussian radial basis functions over the data set with $\sigma=0.005$. In Fig. \ref{henon_pf1} we show the eigenfunction with eigenvalue one of the matrix $\bf P$ capturing the chaotic attractor of Henon map. 

{\underline{\it Van der Pol Oscillator}}: The next step of simulation results is performed with Van der Pol Oscillator.
\begin{eqnarray}
\ddot x=(1-x^2)\dot x-x.
\end{eqnarray}
Time-domain simulation are performed by using discretization time-step of $\Delta t=0.1$ over total time period of $T=10$. The differential equation is solved in MATLAB with $ode45$ solver.
Simulation results from $100$ different randomly chosen initial conditions are generated. For dictionary set we choose $500$ dictionary functions with centers of the dictionary functions determined using $k$-means clustering algorithm with $\sigma=0.1$.

\begin{figure}[h!]
\includegraphics[width=.9\linewidth]{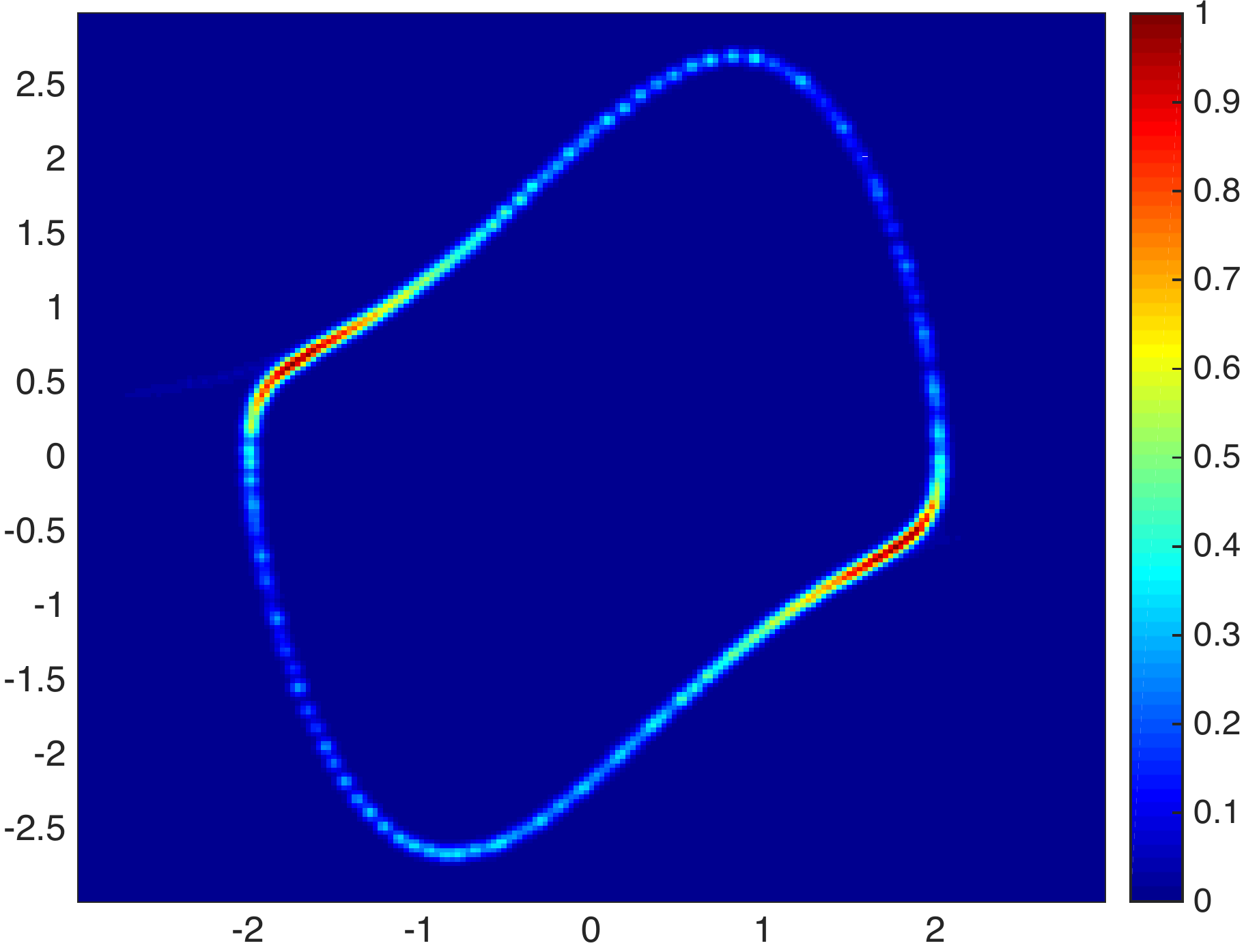}
\centering
\caption{\small {CASE-II: P-F eigenfunction for eigenvalue $1$ for Van der Pol oscillator using NSDMD}}
\label{vanderpol_pf1}
\end{figure}

In Fig. \ref{vanderpol_pf1}, we show the P-F eigenfunctions corresponding to eigenvalue one of the ${\bf P}$ matrix obtained using NSDMD algorithm capturing the limit cycling dynamics of the Vanderpol oscillator. 

\begin{figure}[h!]
\includegraphics[width=.9\linewidth]{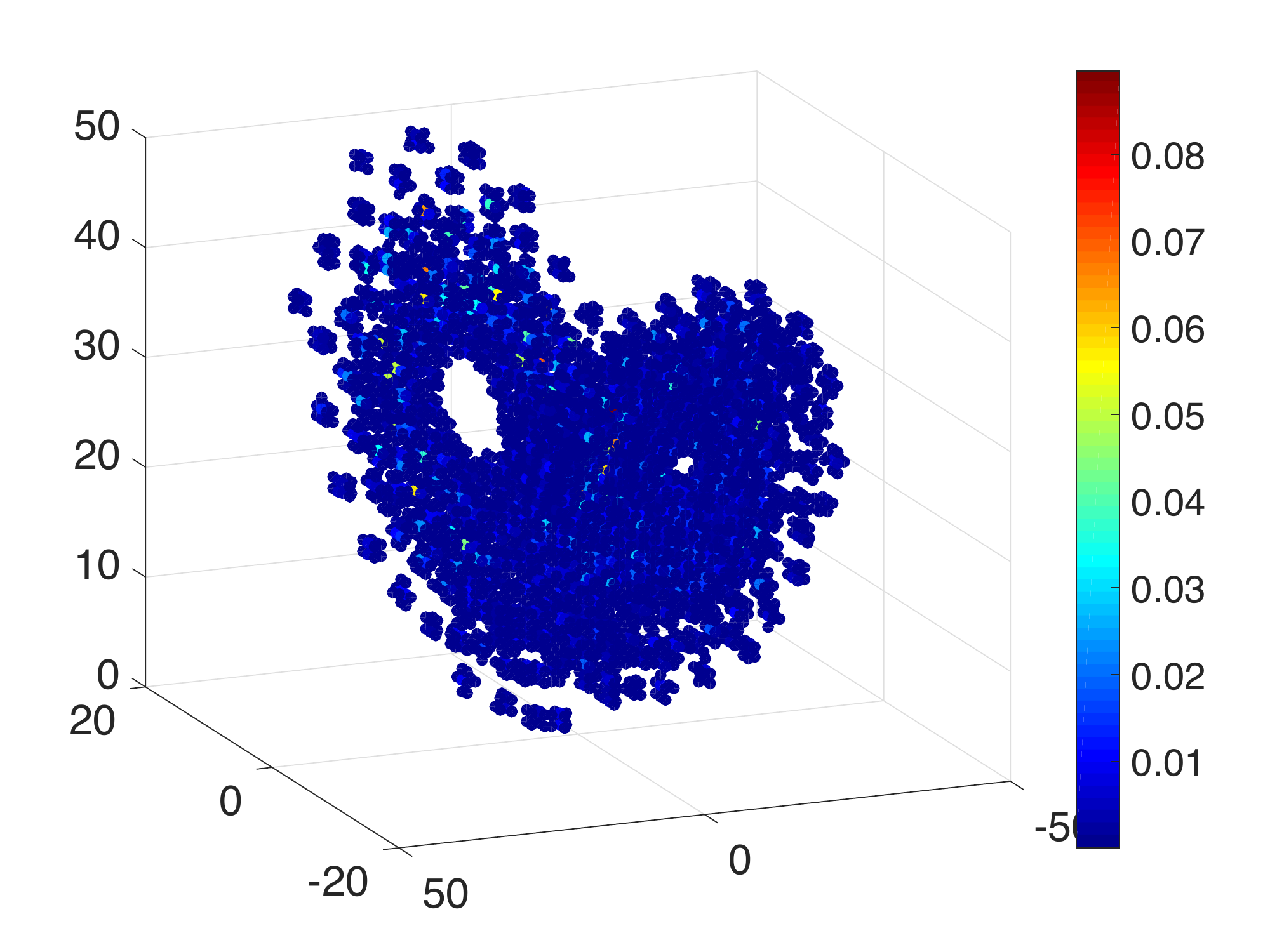}
\centering
\caption{\small {CASE-III: P-F eigenfunction for eigenvalue $1$ for Lorenz attractor using NSDMD}}
\label{lorenz_pf11}
\end{figure}

\begin{figure}[h!]
\includegraphics[width=.9\linewidth]{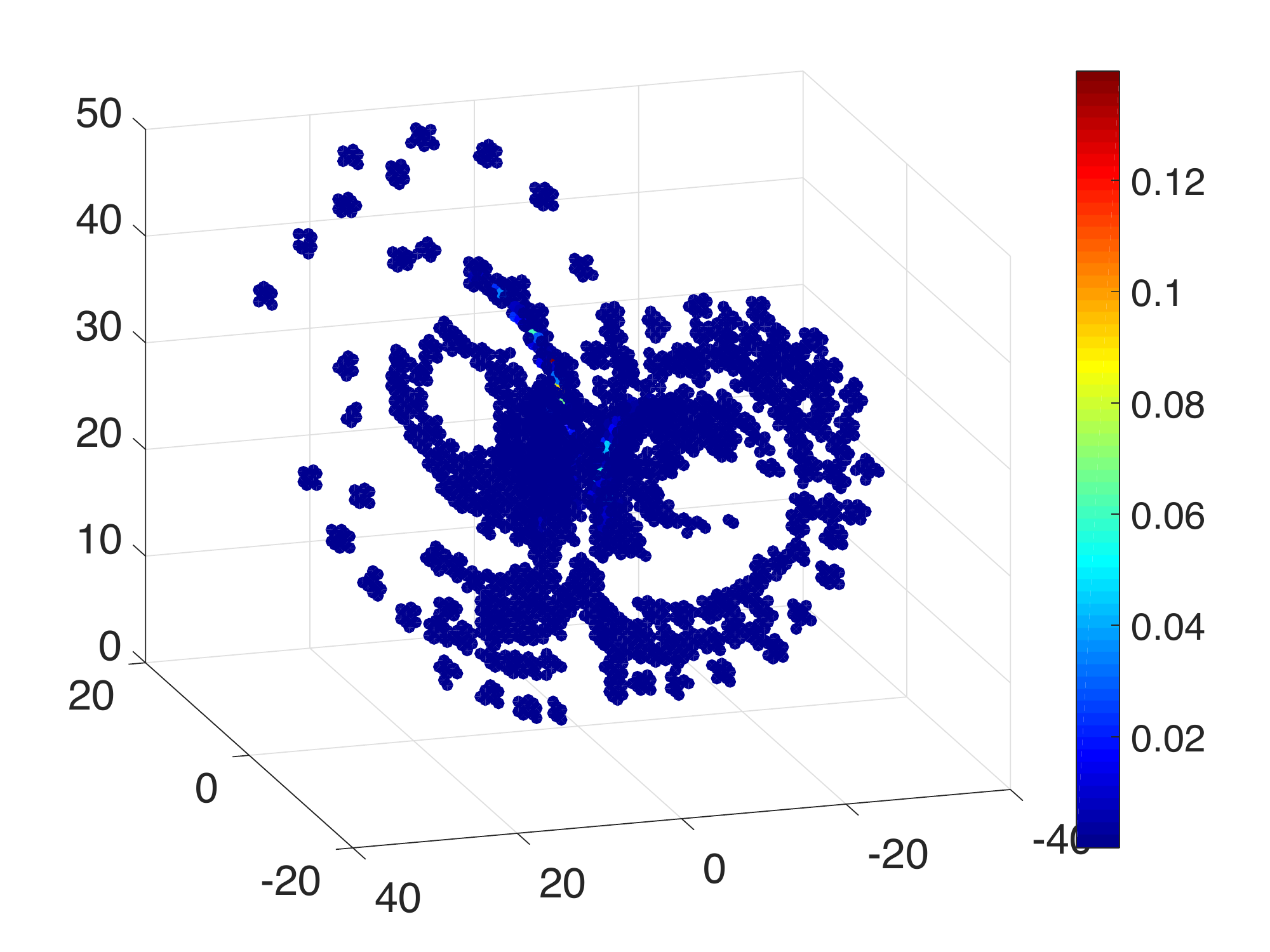}
\centering
\caption{\small {CASE-III: P-F eigenfunction for eigenvalue $0.9762$ for Lorenz attractor using NSDMD}}
\label{lorenz_pf12}
\end{figure}

{\underline{\it Lorenz attractor}}: The simulation results for this example are obtained using optimization formulation in {\bf Case III}.
\begin{eqnarray}
\dot{x} &= a (y - x), \\\nonumber
\dot{y} &= x (b - z) - y, \\\nonumber
\dot{z} &= x y - c z.
\end{eqnarray}
where $a=10$, $b=8/3$ and$c=28$.
Time-domain simulation are performed by using discretization time-step of $\Delta t=0.02$ over total time period of $T=100$. The differential equation is solved in MATLAB with $ode45$ solver.
Simulation results are generated from one initial condition $(1,1,1)$. For dictionary set we choose $500$ dictionary functions with centers of the dictionary functions determined using $k$-means clustering algorithm with $\sigma=0.5$.
In Fig. \ref{lorenz_pf11}, and Fig. \ref{lorenz_pf12}, we show the first two dominant P-F eigenfunctions corresponding to eigenvalue one and $0.9762$ of the ${\bf P}$ matrix obtained using NSDMD algorithm.


\subsection{Stability Certificate: Lyapunov Measure from Data}
\begin{figure}
\includegraphics[width=.8\linewidth]{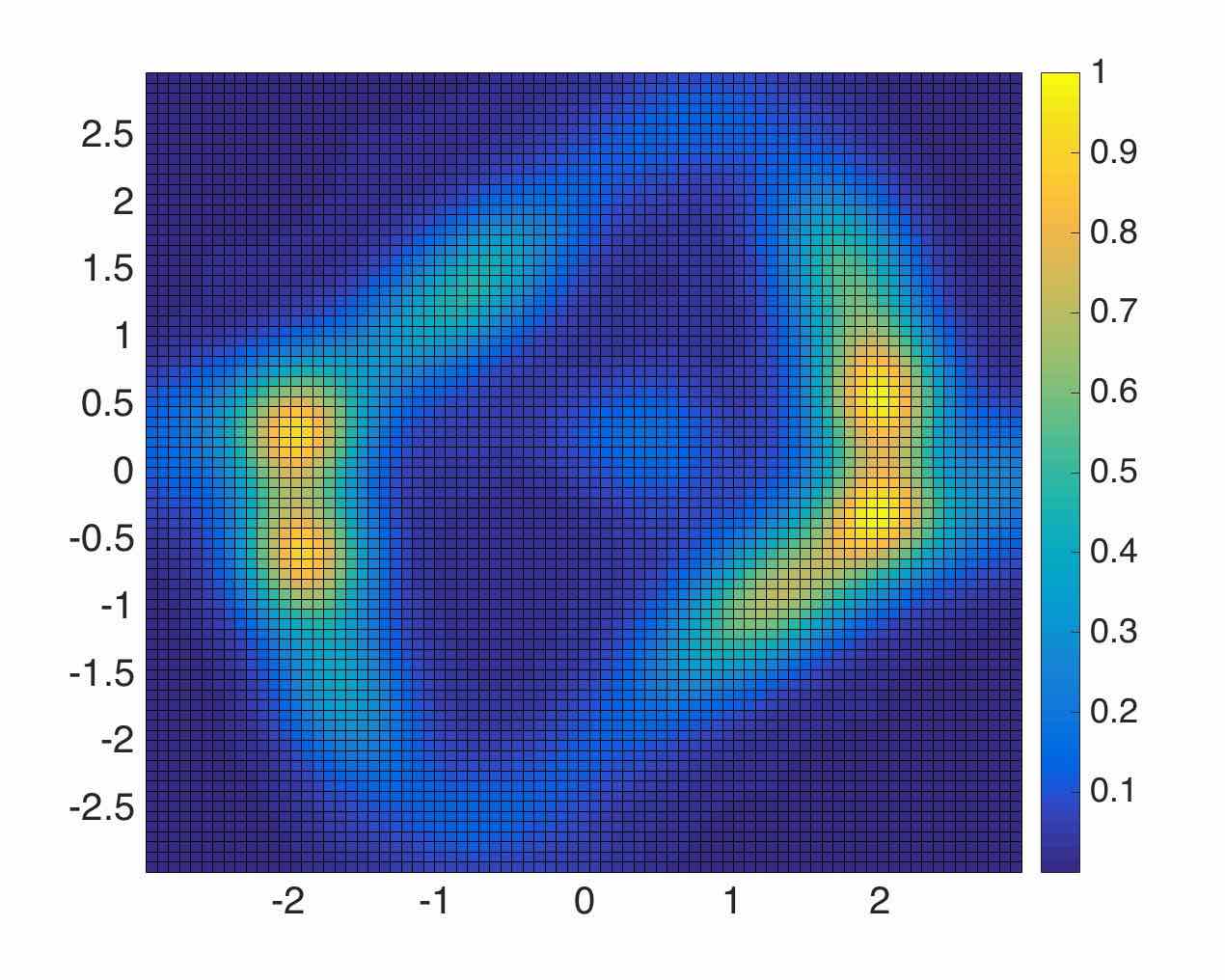}
\centering
\vspace{-0.1in}
\caption{\small{Lyapunov Measure for Van der Pol Oscillator}}
\label{lya}.
\end{figure}

We notice that the finite dimensional Koopman matrix obtained using the DMD or EDMD algorithms are not guaranteed to be stable. For example, in the Van der Pol oscillator, the largest eigenvalue of the $\bf K$ matrix using EDMD is $\lambda=1.001$ and hence unstable. However, the matrix $\bf K$ obtained using NSDMD algorithm is guaranteed to be stable. In fact, the stability certificate in the form of Lyapunov measure can be computed using the procedure outlined in \cite{Vaidya_TAC}.  This stability certificate provides information about the relative amount of time system trajectories spend in the different region of state space before getting absorbed in the attractor set. In Fig. \ref{lya}, we show the plot for the Lyapunov measure for the Van der Pol oscillator example. 

\section{Conclusions}\label{section_conclusion}
We have provided a new algorithm for computing Koopman and P-F eigenfunction from time series data. This proposed algorithm ensure that important properties of the infinite dimensional transfer operators such as positivity and Markov property are preserved in the finite-dimensional approximation. We show via simulation examples that the proposed algorithm can provide a better approximation of the steady-state dynamics regarding eigenfunctions and eigenvalues of the transfer operators. Furthermore, we demonstrate that preserving the positivity property of the finite dimensional approximation is essential to capture the true transient dynamics of the operators.

\bibliographystyle{ieeetr}
\bibliography{ref,ref1,reference}

\end{document}